\documentclass[12pt,twoside]{article}

\date{}

\usepackage{amsmath, amsthm, amssymb,mathrsfs}

\begin{document}
\centerline{\bf International Mathematical Forum, Vol. x, 2008, no. xx, xxx - xxx}

\centerline{}

\centerline{}

\centerline {\Large{\bf Optimal Sequential   Testing of Two
Simple }}

\centerline{}

\centerline{\Large{\bf Hypotheses in Presence of  Control Variables}}

\centerline{}

\centerline{\bf {Andrey Novikov}}

\centerline{}

\centerline{Department of Mathematics}

\centerline{Autonomous Metropolitan University - Iztapalapa}

\centerline{San Rafael Atlixco 186, col. Vicentina}

\centerline{C.P. 09340, Mexico City, Mexico}

\centerline{}

\newtheorem{Theorem}{\quad Theorem}[section]

\newtheorem{Definition}[Theorem]{\quad Definition}

\newtheorem{Corollary}[Theorem]{\quad Corollary}

\newtheorem{Lemma}[Theorem]{\quad Lemma}

\newtheorem{Example}[Theorem]{\quad Example}
\newtheorem{Remark}[Theorem]{\quad Remark}
\begin{abstract} Suppose that at any stage of a statistical experiment  a  control
variable $X$ that affects the distribution of the observed data $Y$
 can be used. The distribution of $Y$ depends on some
unknown parameter $\theta$, and  we consider the classical problem
of testing a simple hypothesis $H_0:\,\theta=\theta_0$ against a
simple alternative $H_1:\,\theta=\theta_1$ allowing the data to be
controlled by $X$, in the following sequential context.

 The experiment starts with assigning a value $X_1$ to the control
variable and observing $Y_1$ as a response. After some analysis, we
choose another value $X_2$ for the control variable, and observe
$Y_2$ as a response, etc.  It is supposed that the experiment
eventually stops, and at that moment a final decision in favour of
$H_0$ or $H_1$ is to be taken.

In this article, our aim  is to characterize the structure of
optimal sequential  procedures, based on this type of data, for
testing a simple hypothesis against a simple alternative.
\end{abstract}



\section{\normalsize  Introduction. Problem Set-Up.}

Let us suppose that at any stage of a statistical experiment  a
"control variable" $X$ that affects the distribution of the observed
data $Y$ can be used. "Statistical" means that  the
distribution of $Y$ depends on some unknown parameter $\theta$, and
we
 have the usual goal of statistical analysis: to obtain some information
about the true value of $\theta$. In this work, we consider the
classical problem of testing a simple hypothesis
$H_0:\,\theta=\theta_0$ versus a simple alternative
$H_1:\,\theta=\theta_1$ allowing the data to be controlled by $X$,
in the following "sequential" context.

 The experiment starts with assigning a value $X_1$ to the control
variable and observing $Y_1$ as a response. After some analysis, we
choose another value $X_2$ for the control variable, and observe
$Y_2$ as a response. Analyzing this, we choose $X_3$ for the third
stage, get $Y_3$, and so on. In this way, we obtain a sequence
$X_1,\dots, X_n$, $Y_1,\dots, Y_n$ of experimental data,
$n=1,2,\dots$. It is supposed that the experiment eventually stops,
and at that moment a final decision in favour of $H_0$ or $H_1$ is to
be taken.

In this article, our aim  is to characterize the structure of
optimal sequential procedures, based on this type of data, for
testing a simple hypothesis against a simple alternative.

Let us write, briefly, $X^{(n)}$ instead of $(X_1,\dots,X_n)$,
$Y^{(n)}$ instead of $(Y_1,\dots,Y_n)$, etc.
 Let us define a (randomized) sequential hypothesis testing
 procedure
 as a triplet $( \chi,\psi, \phi)$ of a a {\em control
policy} $\chi$, a {\em stopping rule} $\psi$, and a {\em decision
rule} $\phi$,
 with $$
\chi=\left(\chi_1, \chi_2, \dots ,\chi_n,\dots\right),
$$
$$\psi=\left(\psi_1, \psi_2, \dots ,\psi_{n},\dots\right),$$
$$\phi=\left(\phi_1, \phi_2, \dots ,\phi_n,\dots\right),$$
where $$\chi_{n}=\chi_{n}(x^{(n-1)},y^{(n-1)})
$$ $n=1,2,\dots$ are supposed to be measurable functions with values
in the space of values of the control variable, and the  functions
$$\psi_n=\psi_n(x^{(n)},y^{(n)}),\quad \phi_n=\phi_n(x^{(n)},y^{(n)})$$
are supposed to be some measurable functions with values in $[0,
1]$.

The interpretation of these functions is as follows.

The experiments starts at stage $n=1$ applying $\chi_1$ to determine
the initial control $x_1$. Using this control, the first data $y_1$
is observed.

At any stage $n\geq 1$:  the value of $\psi_n(x^{(n)},y^{(n)})$ is
interpreted as the conditional probability {\em to stop and proceed
to decision making}, given that that we came to that stage and that
the observations  were $(y_1, y_2, \dots, y_n)$ after the respective
controls $(x_1,x_2,\dots,x_n)$ have been applied. If there is no
stop, the experiments continues to the next stage, defining first
the new control value $x_{n+1}$ by applying the control policy:
$
x_{n+1}=\chi_{n+1}(x^{(n)};y^{(n)})
$
and then taking an additional observation $y_{n+1}$ using control
$x_{n+1}$.

Then the rule $\psi_{n+1}$ is applied to $(x_1,\dots,
x_{n+1};y_1,\dots,y_{n+1})$ in the same way as as above, etc., until
the experiment eventually stops.

It is supposed that when the experiment stops,  a decision {\em to
accept or to reject $H_0$} is to be made. The function
$\phi_n(x^{(n)},y^{(n)})$ is interpreted  as the conditional
probability {\em to reject} the null-hypothesis $H_0$, given that
the experiment stops at stage $n$ being  $(y_1,\dots, y_n)$ the data
vector observed and $(x_1,\dots,x_n)$ the respective controls
applied.

The control policy $\chi$ generates, by the above process, a
sequence of random variables $X_1,X_2,\dots,X_n$, recursively by
$$
X_{n+1}=\chi_{n+1}(X^{(n)},Y^{(n)}).
$$
 The stopping rule $\psi$ generates, by the above
process, a random variable $\tau_\psi$ ({\em stopping time}) whose
distribution is given by
\begin{equation}\label{2a}
P_\theta^\chi(\tau_\psi=n)=E_\theta^\chi (1-\psi_1)(1-\psi_2)\dots
(1-\psi_{n-1})\psi_n.
\end{equation}
Here, and throughout the paper, we interchangeably use  $\psi_n$
both for $$\psi_n(x^{(n)},y^{(n)})$$ and for
$$\psi_n(X^{(n)},Y^{(n)}),$$ and so do we for any other
function  $$F_n=F_n(x^{(n)},y^{(n)}).$$ This does not cause any
problem if we adopt the following agreement: when $F_n$ is under
probability or expectation sign, it is $F_n(X^{(n)},Y^{(n)})$,
otherwise it is $F_n(x^{(n)},y^{(n)})$.


For a sequential testing procedure $(\chi,\psi,\phi)$ let us define
{\em the type I error probability} as
\begin{equation}\label{3a}
\alpha(\chi,\psi,\phi)=P_{\theta_0}(\,\mbox{reject}\,
H_0)=\sum_{n=1}^\infty E_{\theta_0}^\chi(1-\psi_1)
\dots(1-\psi_{n-1})\psi_n\phi_{n}
\end{equation}
and {\em the type II error probability} as
\begin{equation}\label{3c}
\beta(\chi,\psi,\phi)=P_{\theta_1}(\,\mbox{accept}\,
H_0)=\sum_{n=1}^\infty E_{\theta_1}^\chi(1-\psi_1)
\dots(1-\psi_{n-1})\psi_n(1-\phi_{n}).
\end{equation}

Normally, we would like to keep them below some specified levels:
\begin{equation}
\alpha(\chi,\psi,\phi)\leq \alpha \label{1}
\end{equation}
and
\begin{equation}
\beta(\chi,\psi,\phi)\leq \beta \label{2}
\end{equation}
with some $\alpha,\beta\in (0,1)$.

Another important characteristic of a sequential testing procedure
is the {\em average sample number}:
\begin{equation}\label{8aa}
N(\theta;\chi,
\psi)=E_\theta^\chi\tau_\psi=\begin{cases}\sum_{n=1}^\infty
nP_\theta^\chi(\tau_\psi=n), \;\mbox{if}\;
P_\theta^\chi(\tau_\psi<\infty)=1,\cr
\infty\quad\mbox{otherwise}.\end{cases}
\end{equation}

 Our main
goal is minimizing $N(\chi,\psi)=N(\theta_0;\chi,\psi)$ over all
sequential testing procedures subject to (\ref{1}) and (\ref{2}).
Our method is essentially the same that we used in \cite{NovikovSA}
in the problem of sequential testing of two simple hypotheses
without control variables.

In Section \ref{s2}, we reduce the problem of minimizing
$N(\chi,\psi)$ under constraints (\ref{1}) and (\ref{2}) to an
unconstrained minimization problem. The new objective function is
the Lagrange-multiplier function $L(\chi,\psi,\phi)$.

In Section \ref{s3}, we find
$$L(\chi,\psi)=\inf_{\phi}L(\chi,\psi,\phi).$$

In Section \ref{s4}, we minimize $L(\chi,\psi)$ in the class of
truncated stopping rules, i.e. such that $\psi_N\equiv 1$.

In Section \ref{s5}, we characterize the structure of optimal
strategy $(\chi,\psi)$ in the class of non-truncated stopping rules.

In Section \ref{s6}, the likelihood ratio structure for optimal
strategy is given.

In Section \ref{s7}, we apply the results obtained in Section
\ref{s2} -- Section \ref{s5} to minimizing the average sample number
$N(\chi,\psi)$ over all sequential testing procedures subject to
(\ref{1}) and (\ref{2}).

\section{\normalsize Reduction to Non-Constrained
Minimization}\label{s2}

To proceed with minimizing (\ref{8aa}) over the testing procedures
subject to (\ref{1}) and (\ref{2}) let us define the following
Lagrange-multiplier function:
\begin{equation}{\label{4}}
L(\chi,\psi,\phi)=N(\chi,\psi)+\lambda_0\alpha(\chi,\psi,\phi)+\lambda_1\beta(\chi,\psi,\phi)
\end{equation}
where  $\lambda_0\geq 0$ and  $\lambda_1\geq 0$  are some constant
multipliers.

Let $\Delta$ be a class of sequential testing procedures.

The usual relation between the constrained  and the non-constrained
minimization
 is given by the following

\begin{Theorem}\label{t1} Let exist $\lambda_0> 0$ and
$\lambda_1> 0$ and a testing procedure $(\chi^*,\psi^*,\phi^*)\in
\Delta$ such that for any other testing procedure
$(\chi,\psi,\phi)\in \Delta$
\begin{equation}\label{5}
L(\chi^*,\psi^*,\phi^*)\leq L(\chi,\psi,\phi)
 \end{equation}
 holds and such that \begin{equation}\label{6}\alpha(\chi^*,\psi^*,\phi^*)=\alpha
 \quad\mbox{and}\quad
 \beta(\chi^*,\psi^*,\phi^*)=\beta.
 \end{equation}

 Then for any testing procedure $(\chi,\psi,\phi)\in\Delta$ satisfying
 \begin{equation}\label{5bis}
 \alpha(\chi,\psi,\delta)\leq\alpha\quad\mbox{and}\quad
 \beta(\chi,\psi,\delta)\leq\beta
 \end{equation}
 it holds
\begin{equation}\label{5a}
N(\chi^*,\psi^*)\leq  N(\chi,\psi). \end{equation}

 The inequality
in (\ref{5a}) is strict if at least one of the equalities
(\ref{5bis}) is strict.
\end{Theorem}

\begin{proof}
 Let $(\chi,\psi,\phi)\in \Delta$ be any testing procedure
satisfying (\ref{5bis}). Because of (\ref{5}):
$$
L(\chi^*,\psi^*,\phi^*)=N(\chi^*,\psi^*)+\lambda_0\alpha(\chi^*,\psi^*,\phi^*)+\lambda_1\beta(\chi^*,\psi^*,\phi^*)
$$
\begin{equation}\label{5b}\leq
L(\chi,\psi,\phi)=N(\chi,\psi)+\lambda_0\alpha(\chi,\psi,\phi)+\lambda_1\beta(\chi,\psi,\phi)
\end{equation}
\begin{equation}\label{5c}
\leq N(\chi,\psi)+\lambda_0\alpha +\lambda_1\beta,
\end{equation}
where to get the last inequality we used (\ref{1}) and (\ref{2}).

So,
$$
N(\chi^*,\psi^*)+\lambda_0\alpha(\chi^*,\psi^*,\phi^*)+\lambda_1\beta(\chi^*,\psi^*,\phi^*)
\leq N(\chi,\psi)+\lambda_0\alpha +\lambda_1\beta,
$$
and taking into account conditions (\ref{6}) we get from this that
$$N(\chi^*,\psi^*)\leq N(\chi,\psi).
$$

The get the last statement of the theorem we note that if
$N(\chi^*,\psi^*)=N(\chi,\psi)$ then there are equalities in
(\ref{5b})-(\ref{5c}) instead of inequalities which is only possible
if $\alpha(\chi,\psi,\phi)=\alpha$ and
$\beta(\chi,\psi,\phi)=\beta$.
\end{proof}

\section{\normalsize Optimal Decision Rules}\label{s3}

In this section, we start solving the problem of minimizing the
Lagrange-multiplier function $L(\chi,\psi,\phi)$ over all sequential
testing procedures: we first find
$$
\inf_{\phi}L(\chi,\psi,\phi),
$$
and the corresponding decision rule, at which this infimum is
attained.

 Let $I_A$ be the
indicator function of the event $A$.

From this time on, we suppose that for any $n=1,2,\dots,$ the random
variable $Y$, when a control $x$ is applied, has a probability
"density" function
\begin{equation}\label{0} f_\theta(y|x)
 \end{equation}
 (Radon-Nicodym derivative of its distribution) with respect to
 a $\sigma$-finite measure $\mu$ on the respective space.
 We are supposing as well that, at any stage $n\geq 1$, given control values $x_1,x_2,\dots
 x_n$ applied, the observations $Y_1,Y_2,\dots, Y_n$ are
 independent, i.e. their joint probability density function, conditionally on given controls $x_1,x_2,\dots
 x_n$,  can be calculated as
\begin{equation}\label{32}
 f_\theta^n(x_1,\dots,x_n;y_1,\dots, y_n)=\prod_{i=1}^n f_\theta (y_i|x_i),
 \end{equation}
 with respect to the product-measure
 $\mu^n=\mu\otimes\dots\otimes\mu$ of $\mu$ $n$ times by itself.
It is easy to see that any expectation, which uses a control policy
$\chi$, can be expressed as
$$
E_{\theta}^\chi g(Y^{(n)})=\int
g(y^{(n)})f_\theta^{n,\chi}(y^{(n)})d\mu^n(y^{(n)}),
$$
where
$$
f_\theta^{n,\chi}(y^{(n)})=\prod_{i=1}^n f_\theta (y_i|x_i)
$$
with
\begin{equation}\label{31}x_i=\chi_i(x^{(i-1)},y^{(i-1)})\end{equation} for
any $i=1,2,\dots$.

Similarly, for any function $F_n=F_n(x^{(n)},y^{(n)})$ let us define
$$F_n^\chi(y^{(n)})=F_n(x^{(n)},y^{(n)})$$ where $x_1,\dots, x_n$
are defined by (\ref{31}).

 As a first step of minimization of $L(\chi,\psi,\phi)$, let us prove the following
\begin{Theorem} \label{t2} For any $\lambda_0\geq 0$ and $\lambda_1\geq 0$ and for
any sequential testing procedure $(\chi,\psi,\phi)$
\begin{equation}\label{6a}
L(\chi,\psi,\phi)\geq L(\chi,\psi,\phi^*)
\end{equation}
\begin{equation}\label{7}
=N(\chi,\psi)+\sum_{n=1}^\infty \int (1-\psi_1^\chi)\dots
(1-\psi_{n-1}^\chi)\psi_n^\chi
 \min\{\lambda_0f_{\theta_0}^{n,\chi},\lambda_1f_{\theta_1}^{n,\chi}\}
d\mu^n.
\end{equation}
 with
\begin{equation}\label{7a}
\phi^*=(\phi_1^*,\phi_2^*,\dots,\phi_n^*,\dots)
\end{equation}
where
\begin{equation}\label{7aa}
\phi_n^*=I_{\left\{\lambda_0f_{\theta_0}^{n}\leq\lambda_1f_{\theta_1}^{n}\right\}}
\end{equation}
\end{Theorem}
\begin{proof} Inequality (\ref{6a}) is equivalent to
\begin{equation}\label{9}
\lambda_0\alpha(\chi,\psi,\phi)+\lambda_1\beta(\chi,\psi,\phi)\geq\lambda_0\alpha(\chi,\psi,\phi^*)+\lambda_1\beta(\chi,\psi,\phi^*).
\end{equation}

We prove (\ref{9}) by finding a lower bound
for the left-hand side of (\ref{9}) and
proving that this lower bound is
attained at $\phi=\phi^*$ defined by
(\ref{7aa}).

To do this, we will use the following simple

\begin{Lemma}\cite{NovikovSA} \label{l1} Let $\phi, F_1, F_2$ be some measurable functions on
a measurable space with a measure $\mu$, such that
$$
0\leq\phi(x)\leq 1, \quad F_1(x)\geq 0,
\quad F_2(x)\geq 0,\quad
$$
and
$$
\int \min\{F_1(x),F_2(x)\}d\mu(x)<\infty.
$$

Then
$$
\int(\phi(x)F_1(x)+(1-\phi(x))F_2(x))d\mu(x)$$
\begin{equation}\label{7b}
\geq \int \min\{F_1(x),F_2(x)\}d\mu(x)
\end{equation}
with an equality if and only if
\begin{equation}\label{7e}
I_{\{F_1(x)<F_2(x)\}}\leq\phi(x)\leq I_{\{F_1(x)\leq F_2(x)\}}
\end{equation}
$\mu$-almost everywhere.
\end{Lemma}

Starting with the proof of (\ref{9}), let us give to the left-hand
side of it
 the form
 $$
 \lambda_0\alpha(\chi,\psi,\phi)+\lambda_1\beta(\chi,\psi,\phi)
 $$
\begin{equation}\label{7f}
= \sum_{n=1}^\infty \int (1-\psi_1^\chi)\dots
(1-\psi_{n-1}^\chi)\psi_n^\chi
 [\phi_n^\chi\lambda_0f_{\theta_0}^{n,\chi}+(1-\phi_n^\chi)\lambda_1f_{\theta_1}^{n,\chi}]
d\mu^n
 \end{equation}
(see (\ref{3a})).

Applying Lemma 1 to each summand in
(\ref{7f}) we immediately have:
$$
 \lambda_0\alpha(\chi,\psi,\phi)+\lambda_1\beta(\chi,\psi,\phi)
 $$
 \begin{equation}\label{7g}
 \geq \sum_{n=1}^\infty \int (1-\psi_1^\chi)\dots (1-\psi_{n-1}^\chi)\psi_n^\chi
 \min\{\lambda_0f_{\theta_0}^{n,\chi},\lambda_1f_{\theta_1}^{n,\chi}\}
d\mu^n
 \end{equation}
with an equality if
$$
\phi_n=I_{\{\lambda_0f_{\theta_0}^{n}\leq
\lambda_1f_{\theta_1}^{n}\}}=\phi_n^*
$$
for any $n=1,2,\dots$. But in this case the right-hand side of
(\ref{7g}) is
$\lambda_0\alpha(\psi,\phi^*)+\lambda_1\beta(\psi,\phi^*))$, so we
get  (\ref{9}).\end{proof}

\begin{Remark}\label{r2}  It is easy to see, using (\ref{8aa}) and (\ref{7g}), that for any $(\chi,\psi)$ such that 
$
P_{\theta_0}^\chi(\tau_\psi<\infty)
$
the minimum value
$L(\chi,\psi,\phi^*)$ in (\ref{6a}) can be represented as
\begin{equation}\label{10}
L(\chi,\psi)=\sum_{n=1}^\infty \int (1-\psi_1^\chi)\dots
(1-\psi_{n-1}^\chi)\psi_n^\chi\left(nf_{\theta_0}^{n,\chi}+l_n^\chi
\right) d\mu^n,
\end{equation}
where, by definition,
$$l_n=\min\{\lambda_0f_{\theta_0}^n,\lambda_1f_{\theta_1}^n\}.$$
\end{Remark}

Let us denote, for the rest of this article,
$$
s_n^{\psi}=(1-\psi_1)\dots(1-\psi_{n-1})\psi_n\quad\mbox{and}\quad c_n^{\psi}=(1-\psi_1)\dots(1-\psi_{n-1})
$$
for any $n=1,2,\dots$. Respectively,
$$
s_n^{\psi,\chi}=(1-\psi_1^\chi)\dots(1-\psi_{n-1}^\chi)\psi_n^\chi\quad\mbox{and}\quad c_n^{\psi,\chi}=(1-\psi_1^\chi)\dots(1-\psi_{n-1}^\chi)
$$
for any $n=1,2,\dots$.

Let also
$$
C_{n}^{\psi,\chi}=\{y^{(n)}:(1-\psi_1^\chi(y^{(1)}))\dots(1-\psi_{n-1}^\chi(y^{(n-1)}))>0\},
$$
for any $n\geq 2$, and  let $C_1^{\psi,\chi}$ be the space of all  $y^{(1)}$, and finally let
$$
\bar C_{n}^{\psi,\chi}=\{y^{(n)}:(1-\psi_1^\chi(y^{(1)}))\dots(1-\psi_{n}^\chi(y^{(n)}))>0\},
$$
for any $n\geq 1$.

\section{\normalsize Truncated Stopping Rules}\label{s4}

Our next goal is to find a control policy $\chi$ and a stopping rule
$\psi$ minimizing the value of $L(\chi, \psi)$ in (\ref{10}).

In this section, we solve, as an intermediate step,  the problem of
minimization of $L(\chi,\psi)$ over all $\chi$ and $\psi$, where
$\psi\in\Delta^N$,  the class of truncated stopping rules, that is,
\begin{equation}\label{11}
\psi=(\psi_1,\psi_2,\dots,\psi_{N-1},1,\dots).\end{equation}

The following lemma takes over a large part of work of doing this.

\begin{Lemma}\label{l2} Let $r\geq 2$ be any natural number, and let $v_r=v_r(x^{(r)},y^{(r)})$ be
any measurable function. Then
$$
\sum_{n=1}^{r-1}\int
s_n^{\psi,\chi}(nf_\theta^{n,\chi}+l_n^\chi)d\mu^n
+\int
c_r^{\psi,\chi}\left(rf_\theta^{r,\chi}+v_r^\chi\right)d\mu^r
$$
\begin{equation}\label{37}
\geq\sum_{n=1}^{r-2}\int
s_n^{\psi,\chi}(nf_\theta^{n,\chi}+l_n^\chi)d\mu^n
+\int
c_{r-1}^{\psi,\chi}\left((r-1)f_\theta^{r-1,\chi}+v_{r-1}^\chi\right)d\mu^{r-1},
\end{equation}
with
\begin{equation}\label{38}
v_{r-1}=\min\{l_{r-1},f_\theta^{r-1}+R_{r-1}\},
\end{equation}
where
\begin{equation}\label{38b}
R_{r-1}(x^{(r-1)},y^{(r-1)})=\min_{x_r}\int
v_r(x_1,\dots,x_r;y_1,\dots,y_r)d\mu(y_r)
\end{equation}

There is an equality in (\ref{37}) if and only if
\begin{equation}\label{39b}
I_{\{l_{r-1}^\chi<
f_\theta^{r-1,\chi}+R_{r-1}^\chi\}}\leq\psi_{r-1}^\chi\leq
I_{\{l_{r-1}^\chi\leq f_\theta^{r-1,\chi}+R_{r-1}^\chi\}}
\end{equation}
$\mu^{r-1}$-almost everywhere on $C_{r-1}^{\psi,\chi}$, and
\begin{equation}\label{38a}
\int v_r^\chi(y^{(r)})d\mu(y_r)=R_{r-1}^\chi(y^{(r-1)})
\end{equation}
$\mu^{r-1}$-almost everywhere on $\bar C_{r-1}^{\psi,\chi}$.
(We suppose that $R_{r-1}$ defined by (\ref{38b}) is a  measurable function of its
arguments).
\end{Lemma}
\begin{proof} To prove (\ref{37}), it is sufficient to show that
$$
\int
s_{r-1}^{\psi,\chi}((r-1)f_\theta^{{r-1},\chi}+l_{r-1}^\chi)d\mu^{r-1}
+\int
c_r^{\psi,\chi}\left(rf_\theta^{r,\chi}+v_r^\chi\right)d\mu^r
$$
\begin{equation}\label{40}
\geq\int
c_{r-1}^{\psi,\chi}\left((r-1)f_\theta^{r-1,\chi}+v_{r-1}^\chi\right)d\mu^{r-1}.
\end{equation}

By the Fubini theorem, the left-hand side of  (\ref{40}) is equal to
$$
\int
s_{r-1}^{\psi,\chi}((r-1)f_\theta^{r-1,\chi}+l_{r-1}^\chi)d\mu^{r-1}
+\int c_r^{\psi,\chi}\left(\int
\left(rf_\theta^{r,\chi}+v_r^\chi\right)d\mu(y_r)\right)d\mu^{r-1}
$$
\begin{equation}\label{41} =\int
c_{r-1}^{\psi,\chi}[\psi_{r-1}^\chi((r-1)f_\theta^{r-1,\chi}+l_{r-1}^\chi)+
(1-\psi_{r-1}^\chi)\int\left(rf_\theta^{r,\chi}+v_r^\chi\right)d\mu(y_r)]d\mu^{r-1}.
\end{equation}

Because of (\ref{32}),
$$
\int
f_\theta^r(x^{(r)},y^{(r)})d\mu(y_r)=f_\theta^{r-1}(x^{(r-1)},y^{(r-1)}),
$$
 so that the right-hand side of (\ref{41}) transforms to
\begin{equation*}
\int
c_{r-1}^{\psi,\chi}\left[(r-1)f_\theta^{r-1,\chi}+\psi_{r-1}^\chi l_{r-1}^\chi+
(1-\psi_{r-1}^\chi)(f_\theta^{r-1,\chi}+\int v_r^\chi
d\mu(y_r))\right]d\mu^{r-1}
\end{equation*}
\begin{equation}\label{42a}\geq\int
c_{r-1}^{\psi,\chi}\left[(r-1)f_\theta^{r-1,\chi}
+\psi_{r-1}^\chi l_{r-1}^\chi+
(1-\psi_{r-1}^\chi)\left(f_\theta^{r-1,\chi}+R_{r-1}^\chi\right)\right]d\mu^{r-1}
\end{equation}
Applying Lemma \ref{l1} with
\begin{equation*}
\phi=\psi_{r-1}^\chi,\quad
F_1=c_{r-1}^{\psi,\chi}l_{r-1}^\chi,\quad
F_2=c_{r-1}^{\psi,\chi}(f_\theta^{r-1,\chi}+R_{r-1}^\chi),
\end{equation*}
we see that the right-hand side of (\ref{42a}) is greater than or equal to
$$
\int c_{r-1}^{\psi,\chi}\left[(r-1)f_\theta^{r-1,\chi}+\min\{l_{r-1}^\chi,f_\theta^{r-1,\chi}+R_{r-1}^\chi\}\right]d\mu^{r-1}\nonumber
$$
\begin{equation}\label{43}
=\int
c_{r-1}^{\psi,\chi}[(r-1)f_\theta^{r-1,\chi}
+v_{r-1}^\chi]d\mu^{r-1},
\end{equation}
by the definition of $v_{r-1}$ in (\ref{38}).

Moreover, by the same Lemma 1, the right-hand side of (\ref{42a}) is equal to
(\ref{43}) if and only if  $\psi_{r-1}$ satisfies (\ref{39b}) $\mu^{r-1}$-almost everywhere on $C_{r-1}^{\psi,\chi}$.

In addition, there is an equality in (\ref{42a})
if and only if $\chi_{r}$ satisfies (\ref{38a}) $\mu^{r-1}$-almost everywhere on $\bar C_{r-1}^{\psi,\chi}$.
\end{proof}

The following Theorem gives some lower bounds for $L(\chi,\psi)$ when the stopping rule $\psi$ is truncated ($\psi\in\Delta^N$) and characterizes the stopping rules that attain these bounds.

\begin{Theorem}\label{t3} Let $\psi\in\Delta^N$ be any (truncated) stopping rule, and $\chi$ any control policy. Then for any
$1\leq r\leq N-1$ the following inequalities hold true
\begin{equation}\label{46a}
L(\chi,\psi)\geq\sum_{n=1}^{r}\int
s_n^{\psi,\chi}(nf_{\theta_0}^{n,\chi}+l_n^\chi)d\mu^n
+\int
c_{r+1}^{\psi,\chi}\left((r+1)f_{\theta_0}^{r+1,\chi}+V_{r+1}^{N,\chi}\right)d\mu^{r+1}
\end{equation}
\begin{equation}\label{46b}
\geq \sum_{n=1}^{r-1}\int
s_n^{\psi,\chi}(nf_{\theta_0}^{n,\chi}+l_n^\chi)d\mu^n
+\int
c_r^{\psi,\chi}\left(rf_{\theta_0}^{r,\chi}+V_{r}^{N,\chi}\right)d\mu^r,
\end{equation}
where $V_N^N\equiv l_N$, and recursively for $k=N, N-1, \dots 2 $
\begin{equation}\label{48}
V_{k-1}^N=\min\{l_{k-1},f_{\theta_0}^{k-1}+R_{k-1}^N\},
\end{equation}
with
\begin{equation}\label{48a}
R_{k-1}^N=R_{k-1}^N(x^{(k-1)};y^{(k-1)})=\min_{x_{k}}\int V_{k}^N(x_1,\dots,x_k;y_1,\dots,y_k)d\mu(y_{k}).
\end{equation}

 The lower bound in (\ref{46b}) is attained if and only if
\begin{equation}\label{49}
I_{\{l_{k}^\chi< f_{\theta_0}^{k,\chi}+R_k^{N,\chi} \}}\leq\psi_{k}^\chi\leq I_{\{l_{k}^\chi\leq f_{\theta_0}^{k,\chi}+R_k^{N,\chi} \}}
\end{equation}
$\mu^k$-almost everywhere on $C_k^{\psi,\chi}$
and
\begin{equation}\label{49a}
R_k^{N,\chi}(y^{(k)})=\int V_{k+1}^{N,\chi}d\mu(y_{k+1})
\end{equation}
 $\mu^k$-almost everywhere on $\bar C_{k}^{\psi,\chi}$, for any $k=r,\dots, N-1$.
\end{Theorem}
\begin{Remark}\label{r4}
It is supposed in Theorem \ref{t3} and in what follows in this article that all the functions $R_k^N$
defined by (\ref{48a})  are well-defined and
measurable for any $k=1,2,\dots, N$ and for any $N=1,2,\dots$, and that $R_0^N$ defined by (\ref{50b}) below is well defined as well (this is true, for example, if $x_i$ can take only a finite number of values for any $i=1,2,\dots$).
\end{Remark}
\begin{proof}
There is an equality in (\ref{46a}) if $r=N-1$. The rest of the proof immediately follows from Lemma \ref{l2} by induction.
\end{proof}
\begin{Corollary} For any truncated stopping rule $\psi\in\Delta^N$, and for any control rule $\chi$
\begin{equation}\label{49b}
L(\chi,\psi)\geq  1+R_0^N,
\end{equation}
where
\begin{equation}\label{50b}
R_{0}^N=\min_{x_{1}}\int V_{1}^N(x_1;y_1)d\mu(y_{1}).
\end{equation}
The lower bound in (\ref{49b}) is attained if and only if (\ref{49}) is satisfied $\mu^k$-almost everywhere on $C_k^{\psi,\chi}$ and
(\ref{49a}) is satisfied  $\mu^k$-almost everywhere on $\bar C_{k}^{\psi,\chi}$, for any $k=1,2,\dots, N-1$ and,
additionally,
\begin{equation}\label{49e}
R_0^N=\int V_{1}^{N}(\chi_1;y_1)d\mu(y_{1}).
\end{equation}
\end{Corollary}

\begin{Remark}\label{r3} It is obvious that  the testing procedure attaining the lower bound in (\ref{49b})
is optimal among all truncated  testing procedures with
$\psi\in\Delta^N$. But  it only makes practical sense if
$$\min\{\lambda_0,\lambda_1\}> 1 +R_0^N.$$

The reason is  that $\min\{\lambda_0,\lambda_1\}$ can be considered
as "the $L(\chi,\psi)$" function for a trivial sequential test
$(\psi_0,\phi_0)$ which, without taking any observations, makes the
decision $\phi _0=I_{\{\lambda_0\leq\lambda_1\}}$. In this case
there are no  observations ($N(\theta;\psi_0)=0$) and it is easily
seen that
$$L(\psi_0,\phi_0)=\lambda_0\alpha(\psi_0,\phi_0)+\lambda_1\beta(\psi_0,\phi_0)=
\min\{\lambda_0,\lambda_1\}.$$ Thus, the inequality
$$\min\{\lambda_0,\lambda_1\}\leq 1+R_0^N $$
means that the trivial test $(\psi_0,\phi_0)$ is not worse than the
best testing procedure with $\psi$ from $\Delta^N$.

Because of that, we consider
$$V_0^N=\min\{\min\{\lambda_0,\lambda_1\},1+R_0^N\}$$ as the minimum value of $L(\chi,\psi)$ for
$\psi\in\Delta^N$, when taking
no observations is permitted. It is obvious that
this is a particular case of (\ref{48}) with $k=1$, if we define
$l_0\equiv\min\{\lambda_0,\lambda_1\}$ and $f_{\theta_0}^0\equiv 1$.
\end{Remark}

\section{\normalsize Non-Truncated Stopping Rules}\label{s5}

In this section we characterize the structure of general sequential
testing procedures minimizing $L(\chi,\psi)$.

 Let us define for any stopping rule $\psi$ and any control policy
 $\chi$
\begin{equation}\label{50}
L_N(\chi,\psi)=\sum_{n=1}^{N-1}\int
s_n^{\psi,\chi}(nf_{\theta_0}^{n,\chi}+l_n^\chi)d\mu^n+\int c_N^{\psi,\chi}\left(Nf_{\theta_0}^{N,\chi}+l_N^\chi\right)d\mu^{N}.
\end{equation}
This is the Lagrange-multiplier function
corresponding to  $\psi$ truncated  at $N$, i.e. the rule with the
components $\psi^N=(\psi_1,\psi_2,\dots,\psi_{N-1},1,\dots)$:
$$L_N(\chi,\psi)=L(\chi,\psi^N).$$

Because $\psi^N$ is truncated, the results of the preceding section
apply, in particular, the inequalities of Theorem \ref{t3}.

The idea of what follows is to make $N\to\infty$, to obtain some
lower bounds for $L(\chi,\psi)$ from (\ref{46a}) - (\ref{46b}).

And the first question is: what happens to  $L_N(\chi,\psi)$ when
$N\to\infty$?

Let us denote by $\mathscr F$ the set of all strategies
($\chi,\psi$) such that
\begin{equation}\label{50a}
\lim_{n\to\infty} E_{\theta_0}^\chi (1-\psi_1)\dots(1-\psi_{n})=0.
\end{equation}
It is easy to see that (\ref{50a}) is equivalent to
$$
P_{\theta_0}^\chi(\tau_\psi<\infty)=1
$$
(see (\ref{2a})).
\begin{Lemma}\label{l3} For any  strategy $(\chi,\psi)\in \mathscr F$
$$
\lim_{N\to\infty}L_N(\chi,\psi)=L(\chi,\psi).
$$
\end{Lemma}

\begin{proof} Let $L(\chi,\psi)<\infty$, leaving the possibility
$L(\chi,\psi)=\infty$ till the end of the proof. Let us calculate
the difference between $L(\chi,\psi)$ and $L_N(\chi,\psi)$ in order
to show that it goes to zero as $N\to\infty$. By (\ref{50})
$$
L(\psi)-L_N(\psi)=\sum_{n=1}^\infty\int
s_n^{\psi,\chi}(nf_{\theta_0}^{n,\chi}+l_n^\chi)d\mu^n
$$
$$-\sum_{n=1}^{N-1}\int
s_n^{\psi,\chi}(nf_{\theta_0}^{n,\chi}+l_n^\chi)d\mu^n
-\int
c_N^{\psi,\chi}\left(Nf_{\theta_0}^{N,\chi}+l_N^\chi\right)d\mu^n
$$
\begin{equation}\label{51}
=\sum_{n=N}^\infty \int
s_n^{\psi,\chi}(nf_{\theta_0}^{n,\chi}+l_n^\chi)d\mu^n-\int
c_N^{\psi,\chi}\left(Nf_{\theta_0}^{N,\chi}+l_N^\chi\right)d\mu^n.
\end{equation}

The first summand on the right-hand side of (\ref{51}) converges to zero, as $N\to\infty$, being the tail
of a convergent series (this is because $L(\chi,\psi)<\infty$).

We have further
$$\int c_N^{\psi,\chi}l_N^\chi d\mu^n\leq
\lambda_0\int
c_N^{\psi,\chi}f_{\theta_0}^{N,\chi}d\mu^n=\lambda_0E_{\theta_0}^\chi(1-\psi_1)\dots(1-\psi_{N-1})\to 0$$
as $N\to\infty$, because of (\ref{50a}).

It remains to show that
\begin{equation}\label{52}
\int c_N^{\psi,\chi}Nf_{\theta_0}^Nd\mu^n\to 0
\quad\mbox{as}\quad N\to\infty.
\end{equation}

But this is again due to the fact that  $L(\chi,\psi)<\infty$ which
implies that $$E_{\theta_0}^\chi\tau_\psi=\sum_{n=1}^\infty
nP_{\theta_0}(\tau_\psi=n)<\infty.$$ Because this series is
convergent, $\sum_{n=N}^\infty nP_{\theta_0}^\chi(\tau_\psi=n)\to
0$. Thus, using the Chebyshev inequality we have
$$
NP_{\theta_0}^\chi(\tau_\psi\geq N)\leq E_{\theta_0}^\chi\tau_\psi
I_{\{\tau_\psi\geq N\}}=\sum_{n=N}^\infty
nP_{\theta_0}^\chi(\tau_\psi=n)\to 0
$$
as $N\to\infty$, which completes the proof of (\ref{52}).

Let now $L(\chi,\psi)=\infty$.

This means that
$$\sum_{n=1}^{\infty}\int
s_n^{\psi,\chi}(nf_{\theta_0}^{n,\chi}+l_n^\chi)d\mu^n=\infty$$
which immediately implies by (\ref{50}) that
$$L_N(\chi,\psi)\geq\sum_{n=1}^{N-1}\int
s_n^{\psi,\chi}(nf_{\theta_0}^{n,\chi}+l_n^\chi)d\mu^n\to\infty.$$
\end{proof}

The second question is about the behaviour of  the functions $V_r^N$
which participate in the inequalities of  Theorem \ref{t3}, as
$N\to\infty$.

\begin{Lemma}\label{l4} For any $r\geq 1$ and for any  $N\geq r$
\begin{equation}\label{52a}
V_r^N\geq V_r^{N+1}.
\end{equation}
\end{Lemma}

\begin{proof} By induction over $r=N,N-1,\dots,1$.

Let $r=N$. Then by (\ref{48})
$$V_N^{N+1}=\min\{l_N,f_{\theta_0}^N+\min_{x_{N+1}}\int V_{N+1}^{N+1}d\mu(y_{N+1})\}\leq
l_N=V_N^N.
$$

 If we suppose that (\ref{52a}) is satisfied for some $r$, $N\geq
 r>1$, then
$$
V_{r-1}^N=\min\{l_{r-1},f_{\theta_0}^{r-1}+\min_{x_r}\int
V_r^{N}d\mu(y_r)\}$$
$$\geq\min\{l_{r-1},f_{\theta_0}^{r-1}+\min_{x_r}\int V_r^{N+1}d\mu(y_r)\}=V_{r-1}^{N+1}.
$$
Thus, (\ref{52a}) is satisfied for $r-1$ as well, which completes
the induction.
\end{proof}

It follows from Lemma \ref{l4} that for any fixed $r\geq 1$ the
sequence $V_r^N$ is non-increasing. So, there exists
\begin{equation}\label{54}V_r= \lim_{N\to\infty}
V_r^N.
\end{equation}

Now, everything is prepared for passing to the limit, as
$N\to\infty$, in (\ref{46a}) and (\ref{46b}) with $\psi=\psi^N$.

\begin{Theorem}\label{t4} Let $\chi$ be any control policy and $\psi$ any stopping rule. Then for
any $r\geq 1$ the following inequalities hold
\begin{equation}\label{55}L(\chi,\psi)\geq\sum_{n=1}^{r}\int
s_n^{\psi,\chi}(nf_{\theta_0}^{n,\chi}+l_n^\chi)d\mu^n
+\int c_{r+1}^{\psi,\chi}\left((r+1)f_{\theta_0}^{r+1,\chi}+V_{r+1}^\chi\right)d\mu^{r+1}
\end{equation}
\begin{equation}\label{56}
\geq \sum_{n=1}^{r-1}\int
s_n^{\psi,\chi}(nf_{\theta_0}^{n,\chi}+l_n^\chi)d\mu^n+\int
 c_{r}^{\psi,\chi}\left(rf_{\theta_0}^{r,\chi}+V_{r}^\chi\right)d\mu^r,
\end{equation}
where
\begin{equation}\label{56a}
V_r=\min\{l_r,f_{\theta_0}^r+R_r\},
\end{equation}
being
\begin{equation}\label{56b}
R_r=R_r(x^{(r)},y^{(r)})=\min_{x_{r+1}}\int V_{r+1}(x^{(r+1)},y^{(r+1)})d\mu(y_{r+1}).
\end{equation}

In particular, for $r=1$, the following lower bound holds true:
\begin{equation}\label{57a}
L(\chi,\psi)\geq 1+\int V_1(\chi_1,y_1) d\mu(y_1)\geq 1+R_0,
\end{equation}
where, by definition,
$$
R_0=\min_{x_1}\int V_1(x_1,y_1) d\mu(y_1).
$$
 \end{Theorem}
 \begin{proof} Let $(\chi,\psi)\in \mathscr F$ be any strategy.
 Then, by Lemma \ref{l3}, the
left-hand side of (\ref{46a}) tends to $L(\chi,\psi)$ as
$N\to\infty$.

By the Lebesgue monotone convergence theorem, in view of  Lemma
\ref{l4},  passing to the limit  on the right-hand sides of
(\ref{46a}) and  (\ref{46b}) is possible as well. Thus, (\ref{55})
and (\ref{56}) follow.

 Let us now prove  (\ref{56a}), starting from
\begin{equation}\label{57b}
V_{r}^N=\min\{l_{r},f_{\theta_0}^{r}+R_{r}^N\},
\end{equation}
with
\begin{equation}\label{57c}
R_{r}^N=\min_{x_{r+1}}\int V_{r+1}^Nd\mu(y_{r+1})
\end{equation}
(see (\ref{48}) and (\ref{48a}), respectively).

By Lemma \ref{l4}, the left-hand side of (\ref{57b}) tends to $V_r$.
Additionally,
$$
R_r^N=\min_{x_{r+1}}\int V_{r+1}^Nd\mu(y_{r+1})\leq \int
V_{r+1}^Nd\mu(y_{r+1}),
$$
so
\begin{equation*}
\lim_{N\to\infty}R_r^N\leq \lim_{N\to\infty} \int
V_{r+1}^Nd\mu(y_{r+1})= \int V_{r+1}d\mu(y_{r+1})
\end{equation*}
 by the Lebesgue theorem on monotone convergence. Thus,
\begin{equation}\label{57d}
\lim_{N\to\infty}R_r^N\leq \min_{x_{r+1}}\int
V_{r+1}d\mu(y_{r+1})=R_r.
\end{equation}
On the other hand, for any $N\geq 1$,
$$
\int V_{r+1}^Nd\mu(y_{r+1})\geq \int V_{r+1}d\mu(y_{r+1}),
$$
so
$$
R_r^N=\min_{x_{r+1}}\int V_{r+1}^Nd\mu(y_{r+1})\geq
\min_{x_{r+1}}\int V_{r+1}d\mu(y_{r+1})=R_r,
$$
hence
$$
\lim_{N\to\infty} R_r^N\geq R_r.
$$
From this and (\ref{57d}), we get that
$$
\lim_{N\to\infty}R_r^N=R_r.
$$
Therefore, from (\ref{57b}) it follows that
$$
V_r=\lim_{N\to\infty}V_{r}^N= \min\{l_{r},f_{\theta_0}^{r}+R_{r}\},
$$
which proves (\ref{56a}).

 \end{proof}
Let us note now that the right-hand side of (\ref{57a}) coincides
with
$$\inf_{(\chi,\psi)\in\mathscr F}L(\chi,\psi).$$
\begin{Lemma}\label{l5}
\begin{equation}\label{50g}
\inf_{(\chi,\psi)\in\mathscr F}L(\chi,\psi)=1+R_0.
\end{equation}
\end{Lemma}
\begin{proof}
Let us denote
$$
U=\inf_{(\chi,\psi)\in\mathscr F}L(\chi,\psi),\quad U_N=1+R_0^N.
$$

By Theorem 3, for any $N=1,2,\dots$
$$
U_N=\inf_{(\chi,\psi):\psi\in\Delta^N}L(\chi,\psi).
$$
Obviously, $U_N\geq U$ for any $N=1,2,\dots$, so
\begin{equation}\label{50i}
\lim_{N\to\infty}U_N\geq U.
\end{equation}

Let us show first that in fact there is an equality in (\ref{50i}).

 Suppose the contrary, i.e. that $\lim_{N\to\infty}U_N= U+4\epsilon$, with some
 $\epsilon>0$.
 We immediately have from this that
 \begin{equation}\label{50j}U_N\geq U+3\epsilon\end{equation} for all sufficiently large $N$.

On the other hand, by the definition of $U$ there exists a $\psi$
such that $U\leq L(\chi,\psi)\leq U+\epsilon$ and
$(\chi,\psi)\in\mathscr F$.

Because, by Lemma \ref{l3}, $L_N(\chi,\psi)\to L(\chi,\psi)$, as
$N\to\infty$,  we have that
\begin{equation}\label{50k}
L_N(\chi,\psi)\leq U+2\epsilon
\end{equation}
for all sufficiently large $N$ as well. Because, by definition,
$L_N(\chi,\psi)\geq U_N$, we have that
$$
U_N\leq U+2\epsilon
$$
for all sufficiently large $N$, which contradicts (\ref{50j}).

Thus, $$\lim_{N\to\infty}U_N=U.$$

Now, to get (\ref{50g}) we note first that
$$
U=\lim_{N\to\infty}U_N=1+\lim_{N\to\infty}\inf_{x_1}\int
V_1^N(x_1;y_1)d\mu(y_1)$$
$$\leq 1+\inf_{x_1}\int V_1(x_1;y_1)d\mu(y_1)=1+R_0.
$$

On the other hand, by Theorem \ref{t4},
$$
U=\inf_{(\chi,\psi)\in\mathscr F}L(\chi,\psi)\geq 1+R_0,
$$
thus,
$$
U=1+R_0.
$$

\end{proof}

 The following theorem characterizes the structure of
the control- and the stopping-part of optimal sequential testing
procedures.
\begin{Theorem}\label{t5} If there is a strategy $(\chi,\psi)\in \mathscr F$ such that
\begin{equation}\label{60}
L(\chi,\psi)=\inf_{(\chi^\prime,\psi^\prime)\in\mathscr F} L(\chi^\prime,\psi^\prime),
\end{equation}
then
\begin{equation}\label{60b}
I_{\{l_{k}^\chi<
f_{\theta_0}^{k,\chi}+R_{k}^{\chi}\}}\leq\psi_{k}^\chi\leq
I_{\{l_{k}^\chi\leq f_{\theta_0}^{k,\chi}+R_{k}^{\chi}\}}
\end{equation}
$\mu^k$-almost everywhere on $C_k^{\psi,\chi}$, and
\begin{equation}\label{60a}
\int V_{k+1}^{\chi}d\mu(y_{k+1})=R_{k}^{\chi}
\end{equation}
$\mu^k$-almost everywhere on $\bar C_k^{\psi,\chi}$,
for any $k=1, 2\dots$, where $\chi_1$ is defined in such a way that
\begin{equation}\label{60c}
\int V_1^{\chi} d\mu(y_1)=R_{0}.
\end{equation}

On the other hand, if a strategy $(\psi,\chi)$ satisfies (\ref{60b}) $\mu^k$-almost everywhere on $C_k^{\psi,\chi}$,
and satisfies (\ref{60a}) $\mu^k$-almost everywhere on $\bar C_k^{\psi,\chi}$, for any $k=1, 2\dots$, where $\chi_1$ is such that (\ref{60c}) is fulfilled,  then $(\psi,\chi)\in \mathscr F$, and (\ref{60}) holds.
\end{Theorem}

\begin{proof} Let $(\chi,\psi)\in \mathscr F$ be any strategy. By Theorem \ref{t4} for
any fixed $r\geq 1$ the following inequalities hold:
\begin{equation}\label{61bis}
L(\chi,\psi)\geq\sum_{n=1}^{r}\int
s_n^{\psi,\chi}(nf_{\theta_0}^{n,\chi}+l_n^\chi)d\mu^n+\int
c_n^{\psi,\chi}\left((r+1)f_{\theta_0}^{r+1,\chi}+V_{r+1}^\chi\right)d\mu^{r+1}
\end{equation}
\begin{eqnarray}
\label{62bis}
&\geq& \sum_{n=1}^{r-1}\int
s_n^{\psi,\chi}(nf_{\theta_0}^{n,\chi}+l_n^\chi)d\mu^n
+\int
c_r^{\psi,\chi}\left(rf_{\theta_0}^{r,\chi}+V_{r}^\chi\right)d\mu^{r}\\
\nonumber&\geq& \dots\\
\label{63bis}
&\geq& \int\psi_1^\chi(f_{\theta_0}^{1,\chi}+l_1^\chi)d\mu^1+\int
(1-\psi_1^\chi)\left(2f_{\theta_0}^{2,\chi}+V_{2}^\chi\right)d\mu^2\\
\label{64bis}
&\geq& 1+\int V_{1}^\chi d\mu\geq 1+R_0.
\end{eqnarray}

Let us suppose that the right-hand side of (\ref{64bis}) is attained
by some $(\chi,\psi)\in \mathscr F$. This means that there are
equalities in {\em all} of the inequalities (\ref{61bis}) -
(\ref{64bis}). Then, first of all, we get that
\begin{equation}\label{60cc}
R_0=\int V_{1}^\chi(y_1) d\mu(y_1),
\end{equation}
and, successively for  $k=1, 2\dots$, each time applying  Lemma
\ref{l2}, that
\begin{equation}\label{60bb}
I_{\{l_{k}^\chi<
f_{\theta_0}^{k,\chi}+R_{k}^{\chi}\}}\leq\psi_{k}^\chi\leq
I_{\{l_{k}^\chi\leq f_{\theta_0}^{k,\chi}+R_{k}^{\chi}\}}
\end{equation}
$\mu^k$-almost everywhere on $C_k^{\psi,\chi}$, and
\begin{equation}\label{60aa}
\int V_{k+1}^{\chi}(y^{(k+1)})d\mu(y_{k+1})=R_{k}^{\chi},
\end{equation}
$\mu^{k}$-almost everywhere on $\bar C_k^{\psi,\chi}$.
The first part of Theorem \ref{t5} is proved.

To prove the second part, let us suppose that $(\chi,\psi)$
satisfies (\ref{60cc}) - (\ref{60aa}). Applying Lemma \ref{l2}, we see that all the
inequalities in (\ref{62bis})-(\ref{64bis}) are in fact equalities
for $\psi=\psi^r=(\psi_1,\dots,\psi_{r},1,\dots)$.

In particular, this means that there exists
\begin{equation}\label{65bis}
\lim_{r\to\infty}[\sum_{n=1}^{r}\int
s_n^{\psi,\chi}(nf_{\theta_0}^{n,\chi}+l_n^\chi)d\mu^n+\int
c_{r+1}^{\psi,\chi}\left((r+1)f_{\theta_0}^{r+1,\chi}+V_{r+1}^\chi\right)d\mu^{r+1}]=1
+R_0.
\end{equation}
From this, it follows immediately that there exists as well
\begin{equation}\label{66bis}
\lim_{r\to\infty}\sum_{n=1}^{r}\int
s_n^{\psi,\chi}(nf_{\theta_0}^{n,\chi}+l_n^\chi)d\mu^n\leq
1+R_0,
\end{equation}
and that 
\begin{equation}\label{68}
\limsup_{r\to\infty}\int
c_{r}^{\psi,\chi}rf_{\theta_0}^{r,\chi}d\mu^{r}=\limsup_{r\to\infty}(rP_{\theta_0}^\chi(\tau_\psi\geq r))<\infty.
\end{equation}
From (\ref{68}) it follows that $P_{\theta_0}^\chi(\tau_\psi\geq r)\to 0$, as $r\to\infty$, i.e. that $(\chi,\psi)\in \mathscr F$.
Now, the left-hand side of (\ref{66bis}) is $L(\chi,\psi)$ (because
$(\chi,\psi)\in\mathscr F$), and hence
\begin{equation}\label{67bis}
L(\chi,\psi)\leq 1+R_0.
\end{equation}

On the other hand, by virtue of  (\ref{61bis}) - (\ref{64bis})
$
L(\chi,\psi)\geq 1+R_0.
$
From this, and  (\ref{67bis}), we see that
$
L(\chi,\psi)=1+R_0.
$
Because, by Lemma \ref{l5}, $$\inf_{(\chi^\prime,\psi^\prime)\in \mathscr F}
L(\chi^\prime,\psi^\prime)=1+R_0,$$ this proves the second part of Theorem
\ref{t5}.
\end{proof}
\begin{Remark}\label{r9}
Theorem \ref{t5} treats the optimality among strategies which take
at least one observation. If we allow not to take any observation,
there is a possibility that the trivial testing procedure (see
Remark \ref{r3}) gives a better result. It is easy to see that this
happens if
$$
\min\{\lambda_0,\lambda_1\}<1+R_0.
$$
\end{Remark}
\begin{Remark}\label{r5}
In a particular case when the control variable takes only {\em one}
value, $x$, Theorem {\ref{t5}} characterizes the optimal stopping
rule in the problem of testing two simple hypotheses for independent
identically distributed (with density $f_\theta(y|x)$) observations
(see \cite{Lorden}, \cite{NovikovSA}, \cite{NovikovIJPAM}). It is
very well known that the optimal stopping rule is based, in this
particular case, on the likelihood ratio statistic (and  the
resulting test is known as the Sequential Probability Ratio Test
(SPRT) \cite{WaldWolfowitz}). Because of this, we will dedicate the
following section to finding a likelihood structure of the optimal
stopping rule in Theorem \ref{t5}, in the general case of
non-trivial control variables.
\end{Remark}

\section{\normalsize Likelihood Ratio Structure of Optimal
Strategy}\label{s6}

In this section, we will give to the optimal strategy in Theorem
\ref{t5} an equivalent form related to the likelihood ratio process.

Let us start with defining the likelihood ratio:
$$
Z_n=Z_n(x^{(n)},y^{(n)})=\prod_{i=1}^n\frac{f_{\theta_1}(y_i|x_i)}{f_{\theta_0}(y_i|x_i)}.
$$
Let us introduce then the following sequence of functions:
\begin{equation}\label{70a}
\rho_0(z)=g(z)\equiv \min\{\lambda_0,\lambda_1z\},
\end{equation}
and for $k=1, 2, 3, \dots$:
\begin{equation}\label{70}
\rho_k(z)=\min\left\{g(z),1+\min_{x}\int
f_{\theta_0}(y|x)\rho_{k-1}\left(z\frac{f_{\theta_1}(y|x)}{f_{\theta_0}(y|x)}\right)d\mu(y)\right\}
\end{equation}
(we are supposing that all $\rho_k$, $k=0,1,2,\dots$ are
well-defined and measurable functions of $z$).
 It is easy to see that (see (\ref{48}),  (\ref{48a}))
$$
V_N^N=f_{\theta_0}^N\rho_0(Z_N),
$$
and for $k=N-1,N-2,\dots, 1$
\begin{equation}\label{72}
V_k^N=f_{\theta_0}^k\rho_{N-k}(Z_k).
\end{equation}
It is not difficult to see (very much like in Lemma \ref{l4}) that
\begin{equation*}
\rho_k(z)\geq \rho_{k+1}(z)
\end{equation*}
for any $k=0,1,2,\dots$, so there exists
\begin{equation}\label{72a}
\rho(z)=\lim_{n\to\infty}\rho_n(z).
\end{equation}
Using arguments similar to those used in the proof of Theorem
\ref{t4}, it can be shown, starting from (\ref{70}), that
\begin{equation}\label{71}
\rho(z)=\min\left\{g(z),1+R(z)\right\},
\end{equation}
where
\begin{equation}\label{72b}
R(z)=\min_{x}\int
f_{\theta_0}(y|x)\rho\left(z\frac{f_{\theta_1}(y|x)}{f_{\theta_0}(y|x)}\right)d\mu(y).
\end{equation}
Let us pass now to the limit, as $N\to\infty$, in (\ref{72}). We see
that
$$
V_k=f_{\theta_0}^k\rho(Z_k).
$$
Using these expressions in Theorem \ref{t5} we get
\begin{Theorem}\label{t6} If  there exists a strategy $(\chi,\psi)\in \mathscr F$ such that
\begin{equation}\label{160}
L(\chi,\psi)=\inf_{(\chi^\prime,\psi^\prime)\in\mathscr F} L(\chi^\prime,\psi^\prime),
\end{equation}
then
\begin{equation}\label{160b}
I_{\{g(Z_k^\chi)< 1+R(Z_k^\chi)\}}\leq\psi_{k}^\chi\leq
I_{\{g(Z_k^\chi)\leq 1+R(Z_k^\chi)\}}
\end{equation}
$P_{\theta_0}^\chi$-almost sure on \begin{equation}\label{160bb}\{y^{(k)}:\; (1-\psi_1^\chi(y^{(1)}))\dots(1-\psi_{k-1}^\chi(y^{(k-1)}))>0\},\end{equation}
and
\begin{equation}\label{160a}
\int
f_{\theta_0}(y|\chi_{k+1})\rho\left(Z_k^{\chi}\frac{f_{\theta_1}(y|\chi_{k+1})}{f_{\theta_0}(y|\chi_{k+1})}\right)d\mu(y)=R(Z_k^{\chi})
\end{equation}
$P_{\theta_0}^\chi$-almost sure on \begin{equation}\label{160aa}\{y^{(k)}:\; (1-\psi_1^\chi(y^{(1)}))\dots(1-\psi_{k}^\chi(y^{(k)}))>0\},\end{equation} where $\chi_1$ is defined in such a way that
\begin{equation}\label{160c}
\int
f_{\theta_0}(y|\chi_{1})\rho\left(\frac{f_{\theta_1}(y|\chi_{1})}{f_{\theta_0}(y|\chi_{1})}\right)d\mu(y)=R(1).
\end{equation}
On the other hand, if $(\chi,\psi)$ satisfies (\ref{160b}) $P_{\theta_0}^\chi$-almost sure on (\ref{160bb}) and
satisfies (\ref{160a}) $P_{\theta_0}^\chi$-almost sure on (\ref{160aa}), for any $k=1,2,\dots$, where $\chi_1$
satisfies (\ref{160c}),  then $(\chi,\psi)\in \mathscr F$ and $(\chi, \psi)$
satisfies (\ref{160}).
\end{Theorem}
\begin{Remark}\label{r6} It is not difficult to see (very much like in \cite{NovikovIJPAM}) that when \begin{equation}\label{161}1+R(\infty)=1+\lim_{z\to\infty}R(z)>\lambda_0,
\end{equation}
there exist $0<A<B<\infty$ such that $g(z)>1+R(z)$ (see (\ref{160b})) is equivalent to $z\in(A,B)$. By Theorem \ref{t8}, this implies, in particular, that the optimal stopping rule is of an SPRT type: stopping occurs when $Z_n^\chi$ for the first time  exits an interval. Nevertheless, unlike the classical problem of sequential testing, this does not help very much in this case of a statistical experiment with control, because an essential part of the problem is the construction of the optimal control rule (see (\ref{160a})), and there is no apparent way to relate it to the stopping constants $A$ and $B$.

If (\ref{161}) does not hold, the optimal stopping rule is still simple, but may seem somewhat strange. For example, if $1+R(\infty)<\lambda_0$, then the optimal strategy prescribes to stop when, for the first time, $Z_n^\chi$ drops below some $A>0$, and accept $H_0$ at that time. In this case, obviously, the experiment may continue indefinitely, with a large probability, if the alternative hypothesis is true. This does not make much practical sense, and we are not sure that this may ever happen in any testing problem with non-trivial control, but we are unable, generally speaking, to prove that (\ref{161}) is always fulfilled.

The reason why the optimal stopping time may not have a finite expectation under one of the hypotheses lies in the definition of the error probabilities (\ref{3a}) and (\ref{3c}) that do not penalize continuing the experiment indefinitely, and/or in the fact that the average sample number under the alternative hypothesis is not taken into account when minimizing the "risk" (see definition of $L(\chi,\psi)$ in (\ref{4})). Similar phenomenons occur even in the "no-control" case and even when the observations are independent and identically distributed, if the average sample number under one of the hypotheses is disregarded as a criterion of optimization (see \cite{Hawix}). Taking  into account the average sample number under both the null- and the alternative hypothesis remedies this problem (see Remark \ref{r7}) below.
\end{Remark}
\begin{Remark}\label{r7}
Considering as a criterion of optimization, instead of $N(\chi,\psi)=E_{\theta_0}^\chi \tau_\psi$ in (\ref{4}), a weighted sum of the two average sample numbers:
$$
N(\chi,\psi)=\pi_0E_{\theta_0}^\chi \tau_\psi+\pi_1E_{\theta_1}^\chi \tau_\psi,
$$
where $\pi_0$  and $\pi_1$ are some positive numbers, leads to a Bayesian problem of sequential testing in the present context. There are almost evident modifications of Theorems \ref{t3}, \ref{t4}, \ref{t5} and \ref{t6} giving solutions to the respective Bayesian problems as well. For example, instead of (\ref{48}) it should be used
\begin{equation}\label{48mod}
V_{k-1}^N=\min\{l_{k-1},\pi_0f_{\theta_0}^{k-1}+\pi_1f_{\theta_1}^{k-1}+R_{k-1}^N\},
\end{equation}
(\ref{49}) should be modified to
\begin{equation}\label{49mod}
I_{\{l_{k}^\chi< \pi_0f_{\theta_0}^{k,\chi}+\pi_1f_{\theta_1}^{k,\chi}+R_k^{N,\chi} \}}\leq\psi_{k}^\chi\leq I_{\{l_{k}^\chi\leq \pi_0f_{\theta_0}^{k,\chi}+\pi_1f_{\theta_1}^{k,\chi}+R_k^{N,\chi} \}},
\end{equation}
etc., etc.
\end{Remark}

\section{\normalsize Application to the Conditional Problem}\label{s7}

In this section, we apply the results obtained in the preceding
sections to minimizing the average sample size
$N(\chi,\psi)=E_{\theta_0}^\chi \tau_\psi$ over all sequential
testing procedures  with error probabilities not exceeding some
prescribed levels.

Combining Theorems \ref{t1}, \ref{t2} and  \ref{t6}, we immediately
have the following
\begin{Theorem}\label{t8} Let $(\chi,\psi)$ satisfy
the conditions of Theorem \ref{t6}, and let $\phi$ be defined by
\begin{equation}\label{7.1}
   \phi_n=I_{\left\{\lambda_0f_{\theta_0}^{n}\leq\lambda_1f_{\theta_1}^{n}\right\}}
\end{equation}
for $n=1,2,\dots$.

Then for any sequential testing procedure
$(\chi^\prime,\psi^\prime,\phi^\prime)$ such that
\begin{equation}\label{7.2}
    \alpha(\chi^\prime,\psi^\prime,\phi^\prime)\leq
    \alpha(\chi,\psi,\phi)\quad\mbox{and}\quad \beta(\chi^\prime,\psi^\prime,\phi^\prime)\leq
    \beta(\chi,\psi,\phi)
\end{equation}
it holds
\begin{equation}\label{7.3}
    N(\chi^\prime,\psi^\prime)\geq N(\chi,\psi).
\end{equation}

The inequality in (\ref{7.3}) is strict if at least one of the
inequalities in (\ref{7.2}) is strict.

If there are equalities in all of the  inequalities in (\ref{7.2})
and (\ref{7.3}), then $(\chi^\prime,\psi^\prime)$ satisfies the conditions of Theorem \ref{t6} as well (with $\chi^\prime$ instead of $\chi$ and
$\psi^\prime$ instead of   $\psi$).
\end{Theorem}
\begin{proof} The only thing  to be proved is the last
assertion.

Let us suppose that  $$\alpha(\chi^\prime,\psi^\prime,\phi^\prime)=
    \alpha(\chi,\psi,\phi),$$
    $$\beta(\chi^\prime,\psi^\prime,\phi^\prime)=
    \beta(\chi,\psi,\phi),$$ and $$N(\chi^\prime,\psi^\prime)=
    N(\chi,\psi).$$
Then, obviously,
\begin{equation}\label{7.4}
    L(\chi,\psi,\phi)=L(\chi,\psi)=L(\chi^\prime,\psi^\prime,\phi^\prime)\geq
    L(\chi^\prime,\psi^\prime)
\end{equation}
(see (\ref{4})) and Remark \ref{r2}.

By Theorem \ref{t6}, there can not be a strict inequality in the last
inequality in (\ref{7.4}), so
$L(\chi,\psi)=L(\chi^\prime,\psi^\prime)$. From Theorem \ref{t6} it
follows now that $(\chi^\prime,\psi^\prime)$ satisfies (\ref{160b})
-- (\ref{160c}) as well.
\end{proof}

 \vspace{4mm}
 \noindent{\normalsize \bf ACKNOWLEDGMENTS}\vspace{2mm}

 The author greatly appreciates the support of the Autonomous
 Metropolitan University, Mexico City, Mexico, where this work was
 done, and the support of the National System of Investigators (SNI)
 of
 CONACyT, Mexico.

This work is also  supported by Mexico's CONACyT Grant no.
CB-2005-C01-49854-F.\vspace{3mm}


\end{document}